\documentclass[3p,times]{elsarticle}
\usepackage{amsmath}
\usepackage{mathrsfs}
\usepackage{graphicx}
\usepackage{amssymb}
\usepackage[figuresright]{rotating}
\usepackage{multirow}
\usepackage{multicol}
\usepackage{float}
\usepackage{subfigure}
\usepackage{color}
\usepackage{epstopdf}

\newtheorem{theorem}{Theorem}[section]

\newtheorem{remark}{Remark}[section]

\newtheorem{example}{Example}[section]

\begin{document}

\begin{frontmatter}

\title{The GMRES method for solving the large indefinite least squares problem via an accelerated preconditioner}

\author{Jun Li$^{a}$\footnote{Corresponding author. Email: junli026430@163.com.}, Lingsheng Meng$^b$}
\address{$^a$~School of Science, Lanzhou University of Technology, Lanzhou, Gansu, 730050, P. R. China\\
$^b$~College of Mathematics and Statistics, Northwest Normal University, Lanzhou, Gansu, 730070, P. R. China}

\begin{abstract}
In this research, to solve the large indefinite least squares problem, we firstly transform its normal equation into a sparse block three-by-three linear systems, then use GMRES method with an accelerated preconditioner to solve it. The construction idea of the preconditioner comes from the thought of Luo et.al [Luo, WH., Gu, XM., Carpentieri, B., BIT 62, 1983-2004(2022)], and the advantage of this is that the preconditioner is closer to the coefficient matrix of the block three-by-three linear systems when the parameter approachs zero. Theoretically, the iteration method under the preconditioner satisfies the conditional convergence, and all eigenvalues of the preconditioned matrix are real numbers and gathered at point $(1,0)$ as parameter is close to $0$. In the end, numerical results reflect that the theoretical results is correct and the proposed preconditioner is effective by comparing with serval existing preconditioners.
\end{abstract}


\end{frontmatter}

\section{Introduction}
\pagestyle{plain} \setcounter{equation}{ 0}
\renewcommand{\theequation}{1.\arabic{equation}}
In this paper, we consider to solve the indefinite least squares (ILS) problem
\begin{equation} \label{1.1}
\mathrm{ILS:}~\min_{x\in \mathbb{R}^{n}}(b-Ax)^TJ(b-Ax),
\end{equation}
where $A\in \mathbb{R}^{m\times n}$ is a large sparse matrix, $m\geq n$,  $b\in \mathbb{R}^{m}$ and $J=\mbox{diag}(I_p,-I_q)$ is the signature matrix, $p+q=m$. The ILS problem will become to the standard least squares (LS) problem if $p = 0$ or $q = 0$, and its quadratic form is definite; if $pq > 0$, the essence of the problem (\ref{1.1}) is to minimize an indefinite quadratic form with respect to the signature matrix $J$, its normal equation is
\begin{equation}\label{1.2}
A^TJAx=A^TJb.
\end{equation}
Obviously, the Hessian matrix for the ILS problem ({\ref{1.1}}) is $2A^T JA$. And the symmetric and positive definiteness (SPD) of the $A^T JA$ determines the uniqueness of (\ref{1.1}).
In actual applications, such as  the total least squares problem \cite{Golub1980,Huffel1991},  the area of optimization known as $H^\infty$ smoothing \cite{Hassibi1993,Sayed1996} and so on, all of them can be transformed into the ILS problem. This kind of problem was firstly proposed in \cite{Chan1998}, then many solvers start researching it \cite{Bojan2021,Diao2019,Li2018,Liu2011,Liu2013,Liu2014,Xu2004}. For the small and dense ILS problem (\ref{1.1}), authors in \cite{Bojan2003,Bojan2021,Chan1998,Xu2004} considered adopting the direct method to solve it, whose thought is mainly based on the $\rm QR$ decomposition by using $\rm Cholesky$ factorization or $\rm Householder$ transformation. However, iteration method may be more attractive for large and sparse ILS problem. In 2011, Liu et.al \cite{Liu2011} constructed a preconditioner based on a partitioning of $A$ for sparse ILS problem and then used the preconditioned conjugate gradient methods for solving the ILS problem (\ref{1.1}). Thereafter, they also established the incomplete hyperbolic Gram-Schmidt-based preconditioners \cite{Liu2013} by using the  incomplete hyperbolic classical Gram–Schmidt or incomplete hyperbolic modified Gram–Schmidt method to approximate the matrix $A^TJA$. In 2014, a block SOR iteration method with relaxation parameter was proposed to solve the ILS problem \cite{Liu2014}. In order to further improve the convergence of SOR iteration method, Song et.al \cite{Song2020} proposed the USSOR iteration method for solving the block $3\times3$ linear systems derived from the ILS problem in 2020. Zhang and Li \cite{Zhang2023} in 2023 gave a splitting-based randomized iteration methods to solve the ILS problem.
In 2025, Meng et.al \cite{Meng2025} proposed a variable Uzawa method to solve the ILS problem. These literatures are too numerous to be listed, impossible to enumerate.

Among iteration methods, Krylov subspace methods (such as GMRES) have attracted wide attention because of its good convergence. However, its performance can be better showed with a suitable preconditioner for solving sparse linear systems.  In 2025, Xin et.al \cite{Xin2025} established three block preconditioners for solving the block three-by-three linear systems derived from the ILS problem (\ref{1.1}). They first partitioned $A$ and $b$ in the ILS problem (\ref{1.1}) as
$$
A=\left(
\begin{array}{c}
A_1\\
A_2
\end{array}
\right),~~b=\left(
\begin{array}{c}
b_1\\
b_2
\end{array}
\right),
$$
where $A_1 \in \mathbb{R}^{p\times n}$ is of full column rank, $A_2 \in \mathbb{R}^{q\times n}$, $b_1 \in \mathbb{R}^{p}$ and $b_2 \in \mathbb{R}^{q}$. Because  $A^TJA$ is SPD, $A_1^TA_1-A_2^TA_2\succ0$. Next, the normal equation (\ref{1.2}) can be  written as the following linear systems:
\begin{equation}\label{1.3}
\mathcal{A}\mathbf{u}=\left(
\begin{array}{ccc}
I & A_1 & 0 \\
A_1^T & 0 & -A_2^T\\
0  & A_2 & I
\end{array}
\right)\left(\begin{array}{c}
\delta_1 \\
x \\
\delta_2
\end{array}
\right)=\left(\begin{array}{c}
b_1 \\
0\\
b_2
\end{array}
\right)=\mathbf{b},
\end{equation}
where $\delta=(\delta_1;\delta_2)=b-Ax$.
Furthermore, to ensure the non-singularity of the diagonal block, they consider the forms of the above linear systems:
\begin{eqnarray}\label{1.4}
\mathcal{\hat{A}}\mathbf{u}=\left(
\begin{array}{ccc}
I &A_1 & 0 \\
0& A_1^T A_1 & A_2^T\\
0  & A_2 & I
\end{array}
\right)\left(\begin{array}{c}
\delta_1 \\
x\\
\delta_2
\end{array}
\right)=\left(\begin{array}{c}
b_1 \\
A_1^Tb_1\\
b_2
\end{array}
\right)=\mathbf{\hat{b}}.
\end{eqnarray}
It is obvious that the above block linear systems (\ref{1.4}) is large and sparse and its coefficient matrix is nonsymmetric indefinite. Then, three block preconditioners was proposed:
\begin{eqnarray}\label{1.5}
{BS}_1=\left(
\begin{array}{ccc}
I &0 & 0 \\
0& A_1^T A_1 & 0\\
0  & 0 & I
\end{array}
\right),\quad
{BS}_2=\left(
\begin{array}{ccc}
I &0 & 0 \\
0& A_1^T A_1 & A_2^T\\
0  & 0 & I
\end{array}
\right),\quad
{BS}_3=\left(
\begin{array}{ccc}
I &A_1 & 0 \\
0& A_1^T A_1 & 0\\
0  & 0 & I
\end{array}
\right).
\end{eqnarray}
Next, the authors analyzed the spectral properties of the preconditioned matrices and verified the efficiency of the preconditioners by numerical experiments; in addition, literature also indicates that $BS_2$ preconditioner is the most effective of the three preconditioners. Thereafter, Li et.al \cite{Li2024} established a block upper triangular (BUT) preconditioner for solving the linear systems (\ref{1.4}), where the preconditioner $P_{BUT}$ is defined as follows:
\begin{align}\label{BUT}
P_{BUT}=\left(
\begin{array}{ccc}
I &A_1 & 0 \\
0& A_1^T A_1 & A_2^T\\
0  & 0 & I
\end{array}
\right).
\end{align}
Similarly, authors also proved the efficiency of the BUT preconditioner. Besides, no existing literature has addressed ILS problem (\ref{1.3}) by converting it into block linear systems and applying preconditioned Krylov subspace methods based on the coefficient matrix structure for solving ILS problem. Therefore, we will study and give an efficient preconditioning method for solving block linear systems derived from ILS problem.

In 2022,  to solve saddle point problem, Luo et.al \cite{Luo2022} considered a dimension expanded preconditioning technique. The main idea is as follows: firstly, the saddle point problem
$$
\left(
\begin{array}{cc}
\hat{A} &\hat{B} \\
\hat{C}& 0
\end{array}
\right)
\left(
\begin{array}{c}
\hat{x}_1 \\
\hat{x}_2
\end{array}
\right)=
\left(
\begin{array}{c}
\hat{b}_1 \\
\hat{b}_2
\end{array}
\right)
$$
is equivalently converted into the systems:
\begin{eqnarray}\label{1.6}
\left(
\begin{array}{ccc}
I&0&I\\
\hat{B}&\hat{A} & 0\\
I&\hat{C}& I
\end{array}
\right)
\left(
\begin{array}{c}
\hat{x}_2 \\
\hat{x}_1\\
\hat{x}_3
\end{array}
\right)=
\left(
\begin{array}{c}
0\\
\hat{b}_1 \\
\hat{b}_2
\end{array}
\right),
\end{eqnarray}
where $\hat{x}_2=-\hat{x}_3$; next, a block preconditioner was proposed for solving the above the $3\times 3$ block linear systems:
\begin{eqnarray}\label{1.7}
\mathcal{P}=\left(
\begin{array}{ccc}
I&0&\alpha I\\
\hat{B}&\hat{A} & 0\\
I&\hat{C}& I
\end{array}
\right),
\end{eqnarray}
here $\alpha $ is a parameter close to $1$; furthermore, theoretical analysis shows that the all eigenvalues of the preconditioned matrix converge to $1$ as $\alpha\rightarrow 1$; in the end, numerical experiments verify the favourable convergence properties of the block preconditioner. In fact, the parameter is closer to 1, the more clustering  the spectrum of the preconditioned matrix is, which also ensures the good convergence of the preconditioning method.

To put it simply, the above thought is about the equivalent transformation of the solved initial linear systems, and then the preconditioner (\ref{1.7}) is established by modifying a sub-block of coefficient matrix in (\ref{1.6}). Inspired by the thought in \cite{Luo2022}, in this paper, instead of directly solving ILS problem (\ref{1.1}), we solve the three-by-three block linear systems (\ref{1.3}). Similarly, we establish a block preconditioner close to coefficient matrix of (\ref{1.3}). When extra parameter $\alpha\rightarrow0$, the preconditioned matrix has clustering eigenvalues. Finally, numerical experiments also verify the efficiency of the proposed preconditioner.

The paper is organized as follows.  a block preconditioner is established for solving the three-by-three block linear systems (\ref{1.3}) in Section 2. Section 3 gives some theoretical result, among which are the convergence of the iteration method under the preconditioner, the eigenvalues distribution of the preconditioned matrix. Numerical experiments in Section 4 show that the proposed preconditioner is effective. Section 5 contains the conclusions.

In the end, we introduce some used notations in the paper. Matlab notation $[x;y;z]$ represents the column vector $[x^T,y^T,z^T]^T$. For a given nonsingular matrix $Z$, $\lambda_{min}(Z)$ and $\rho(Z)$ denotes the minimum eigenvalue and spectral radius of matrix $Z$, respectively, $Z\succ 0$ means that $Z$ is symmetric positive definite (SPD) matrix. The conjugate transpose of a vector $\ell$ is expressed as $\ell^*$.

\section{A new block splitting method}
\pagestyle{plain} \setcounter{equation}{ 0}
\renewcommand{\theequation}{2.\arabic{equation}}

In this section, we establish a block preconditioner, which can be within GMRES method to solve the linear systems (\ref{1.3}).

Combing with idea in \cite{Luo2022}, we give a new splitting of the coefficient matrix $\mathcal{A}$ by modifying a sub-block of coefficient matrix in (\ref{1.3}):
\begin{eqnarray}\label{2.1}
\mathcal{A}=\mathcal{P}-\mathcal{Q}:=\left(
\begin{array}{ccc}
I &A_1 & 0 \\
A_1^T& \alpha I & -A_2^T\\
0  & A_2 & I
\end{array}
\right)-
\left(
\begin{array}{ccc}
0 &0 & 0 \\
0& \alpha I & 0\\
0  & 0 & 0
\end{array}
\right),
\end{eqnarray}
where parameter $\alpha >0$. In order to ensure the nonsingularly of $\mathcal{P}$, $\alpha$ should be close to $0$ such that
$$
\left(\alpha I+A_2^TA_2-A_1^TA_1\right)^{-1}
$$
exists. According to the above matrix splitting, we can establish a new iteration method.\\
{\bf The iteration method.} Assume that the solution of the ILS problem (\ref{1.1}) exists and satisfies uniqueness, its solution can be obtained by solving the block linear systems (\ref{1.3}).   Given an initial guess $\mathbf{u}^{(0)}=(\delta_1^{(0)};x^{(0)}; \delta_2^{(0)}) \in \mathbb{R}^{n+p+q}$, where $(\delta_1^{(0)};\delta_2^{(0)})=b-Ax^{(0)}$, for $k=0,1,2,\ldots $ until iteration sequence $\mathbf{u}^{(k)}=(\delta_1^{(k)};x^{(k)};\delta_2^{(k)}) \in \mathbb{R}^{n+p+q}$ converges, compute
\begin{equation}\label{2.2}
\mathcal{P}\mathbf{u}^{(k+1)}=\mathcal{Q}\mathbf{u}^{(k)}+\mathbf{b}.
\end{equation}

The iteration matrix of the iteration method is
\begin{equation}\label{2.3}
\mathcal{T}=\mathcal{P}^{-1}\mathcal{Q}.
\end{equation}
Then, the condition $\rho(\mathcal{T}) < 1 $ guarantees the convergence of the iteration method.

\begin{remark}
For four preconditioners in (\ref{1.5}) and (\ref{BUT}), $\mathcal{\hat{A}}$ produces the matrix splitting as follows:
$$
\mathcal{\hat{A}}={BS}_1-\mathcal{N}_1={BS}_2-\mathcal{N}_2={BS}_3-\mathcal{N}_3=P_{BUT}-Q_{BUT},
$$
where
$$
\mathcal{N}_1=\left(
\begin{array}{ccc}
0 &-A_1 & 0 \\
0& 0 & -A_2^T\\
0  & -A_2 & 0
\end{array}
\right),\quad \mathcal{N}_2=\left(
\begin{array}{ccc}
0 &-A_1 & 0 \\
0& 0 & 0\\
0  & -A_2 & 0
\end{array}
\right),\quad \mathcal{N}_3=\left(
\begin{array}{ccc}
0 &0 & 0 \\
0& 0 & -A_2^T\\
0  & -A_2 & 0
\end{array}
\right),\quad Q_{BUT}=
\left(
\begin{array}{ccc}
0 &0 & 0 \\
0& 0 &0\\
0  & -A_2 & 0
\end{array}
\right).
$$
It is  not hard to find that all matrices $\mathcal{N}_1,~\mathcal{N}_2$ and $\mathcal{N}_3$ contain two or more nonzero sub-matrices, however, $Q_{BUT}$ and $\mathcal{Q}$  have only one nonzero sub-matrix. Besides, the preconditioner $\mathcal{P}$ has a parameter, and $\mathcal{P}\rightarrow \mathcal{A}$ as $\alpha \rightarrow 0_+$, which means that the preconditioner $\mathcal{P}$ is closer to the coefficient matrix than the other four preconditioners  in (\ref{1.5}) and (\ref{BUT}).
This is expected to make the preconditioned GMRES method under the preconditioner $\mathcal{P}$ more effective than the other four preconditioning methods for solving three-by-three block linear systems.
\end{remark}
%For implementation of the new preconditioner ${P}_{BUT}$, we need to solve sequences of generalized residual equation of the form
%where $[z_1;z_2;z_3]$ and $[r_1;r_2;r_3]$ are the current and the generalized residual vectors, respectively. After some calculations, the procedure for solving the generalized residual equation (\ref{2.6}) are derived as
%\begin{itemize}
%\item[(i)] compute $z_3=r_3$;
%\item[(ii)] solve $(A_1^TA_1)z_2=r_2-A_2^Tz_3$;
%\item[(iii)] compute $z_1=r_1-A_1z_2$.
%\end{itemize}

%In the above algorithm, because the matrix $A_1^TA_1$ is SPD, the Cholesky factorization can be adopted to solve the systems (ii). Besides, the coefficient matrix of the solved linear subsystems after preconditioning is $A_1^TA_1$, but the coefficient matrix of the linear subsystem is $A_1^TA_1-A_2^TA_2$ when the initial systems (\ref{2.2}) is directly  solved, which is expected that the preconditioning method can be solved quickly.

\section{Theoretical  results}
\pagestyle{plain} \setcounter{equation}{ 0}
\renewcommand{\theequation}{3.\arabic{equation}}

In this section, we mainly analyze the convergence of the iteration method and eigenvalues distribution of the preconditioned matrix $\mathcal{P}^{-1}\mathcal{A}$.

%Firstly, we give a lemma, which is useful for proving the theorem.
%\begin{lemma}\cite{6}\label{lem3.1}
%Let the real matrices $A \succ 0$ and $B \succeq 0$. Then $A \succeq B$ if and only if $\rho(A^{-1}B)\leq1$,
%and $A \succ B$ if and only if $\rho(A^{-1}B)<1$, $A\succ B(A \succeq B)$ means $A - B \succ 0(A - B \succeq 0)$.
%\end{lemma}

\begin{theorem}\label{thm3.1}
Solving the block linear systems (\ref{1.3}) to obtain the solution of the ILS problem (\ref{1.1}). The new iteration method denotes in (\ref{2.2}). Then the iteration method is convergent if
$$0< \alpha < \frac{1}{2}\lambda_{min}(S),$$
where $S=A_1^{T}A_1-A_2^{T}A_2 \succ 0$.
\end{theorem}

\noindent{\bf Proof.} Let $(\lambda,\ell)$ with $\ell=(\ell_1;\ell_2;\ell_3)$ be an eigenpair of $\mathcal{P}^{-1}\mathcal{Q}$.  Then we have $\mathcal{P}^{-1}\mathcal{Q}\ell=\lambda \ell$, i.e.,
\begin{eqnarray}\label{3.1}
\left(\begin{array}{ccc}
               0 & 0 & 0 \\
               0 & \alpha I & 0\\
               0 & 0 & 0
             \end{array}\right)
\left(\begin{array}{c}
               \ell_1\\
               \ell_2\\
               \ell_3
             \end{array}\right)=\lambda
\left(\begin{array}{ccc}
               I & A_1& 0 \\
               A_1^T & \alpha I & -A_2^T\\
               0 & A_2 & I
             \end{array}\right)\left(\begin{array}{c}
               \ell_1\\
               \ell_2\\
               \ell_3
             \end{array}\right),
\end{eqnarray}
namely,
\begin{equation}\label{3.2}
\left\{
  \begin{array}{ll}
    \lambda (\ell_1+A_1 \ell_2)=0, \\
    \lambda A_1^T\ell_1+\alpha(\lambda-1)\ell_2-\lambda A_2^T\ell_3=0, \\
    \lambda(A_2\ell_2+\ell_3)=0.
  \end{array}
\right.
\end{equation}
If $\lambda=0$, it is nothing to prove. Next, we consider the case of $\lambda \neq 0$. In this case, we can obtain by (\ref{3.2}) that  $\ell_1=-A_1 \ell_2$ and $\ell_3 =-A_2\ell_2$, submitting them into the second equation of (\ref{3.2}) gives
\begin{equation}\label{3.3}
\lambda(A_1^T A_1-A_2^T A_2) \ell_2+\alpha(1-\lambda)\ell_2=0.
\end{equation}
Assume that $\lambda=1$, i.e., $(A_1^T A_1 -A_2^T A_2) \ell_2=0$, which combing $A_1^TA_1-A_2^TA_2\succ0$ means that $\ell_2=0$, furthermore, we have $\ell_1=0$, $\ell_3=0$, which is a contradiction. Next, we consider $\lambda \neq 1.$ It is obvious that $\ell_2 \neq 0,$ multiplying (\ref{3.3}) from left by $\frac{\ell_2^*}{\ell_2^*\ell_2}$ and letting
\begin{align}\label{3.4}
\tilde{a}:=\frac{\ell_2^*\left(A_1^T A_1-A_2^T A_2\right)\ell_2}{\ell_2^*\ell_2}
\end{align}
obtain
\begin{align}\label{3.5}
(\alpha-\tilde{a})\lambda=\alpha,
\end{align}
where $\tilde{a}>0$ (because $A_1^TA_1-A_2^T A_2$ is SPD). Furthermore, we derive the expression of $\lambda$:
$$
\lambda=\frac{\alpha}{\alpha-\tilde{a}}.
$$

In order to ensure the convergence of the iteration method, the condition $\mid \lambda \mid<1$ needs to be satisfied. After simple calculation, if $0< \alpha < \frac{1}{2}\lambda_{min}(A_1^TA_1-A_2^T A_2)$, $\mid \lambda \mid<1$, i.e., the iteration method is convergent. $\hfill\Box$

\begin{remark}
According to theorem \ref{thm3.1}, we know that the iteration method satisfies the convergence as $0< \alpha < \frac{1}{2}\lambda_{min}(A_1^TA_1-A_2^T A_2)$. In fact, section 2 shows that $\mathcal{P}$ is more close to coefficient matrix $\mathcal{A}$ when $\alpha \rightarrow 0_+$. Based on these truths, we can select parameter $\alpha$ close to $0$ such that the iteration method is always convergent.
\end{remark}

Next, we  will discuss the eigenvalues distribution of the preconditioned matrix $\mathcal{P}^{-1}\mathcal{A}$.
\begin{theorem}\label{thm3.2}
Assume that the conditions of theorem \ref{thm3.1} are satisfied. The new preconditioner is denoted in (\ref{2.1}). Then all eigenvalues of the preconditioned matrix $\mathcal{P}^{-1}\mathcal{A}$ are real numbers and  gather at point $(1,0)$ as $\alpha \rightarrow 0_+$.
\end{theorem}
\noindent{\bf Proof.} Let $\theta$ be a eigenvalue of the preconditioned matrix $\mathcal{P}^{-1}\mathcal{A}$. Because $\mathcal{P}^{-1}\mathcal{A}=\mathcal{P}^{-1}(\mathcal{P}-\mathcal{Q})=I-\mathcal{T}$, it follows that $\theta=1-\lambda$, i.e.,
$$
\theta=\frac{\tilde{a}}{\tilde{a}-\alpha}=1+\frac{\alpha}{\tilde{a}-\alpha},
$$
where $\tilde{a}$ denotes in (\ref{3.4}), which reveals that $\tilde{a}$ is real number, i.e., $\theta$ is real eigenvalue. Besides, when $\alpha\rightarrow 0_+$, $\tilde{a}-\alpha>0,$ therefore, $\theta>0$ and $\theta\rightarrow 1$. This completes the proof. \qed

\begin{remark}
Theorem \ref{thm3.2} shows that the eigenvalues of the preconditioned matrix $\mathcal{P}^{-1}\mathcal{A}$ are gathered at point $(1,0)$ when $\alpha\rightarrow 0_+$; theorem \ref{thm3.1} reveals that if $0< \alpha < \frac{1}{2}\lambda_{min}(A_1^TA_1-A_2^T A_2)$, the iteration method under the preconditioner is convergent. Combing these truths, the parameter $\alpha$ should be selected close to $0$, which will leads to good properties of the preconditioner. This also gives the direction of parameter selection in theory. This perspective will be substantiated through computational testing.
\end{remark}

\begin{remark}
When the sparse linear systems $\mathcal{A}\mathbf{u}=\mathbf{b}$ in (\ref{1.3}) is solved using Krylov subspace methods, the inherent sparsity of the coefficient matrix $\mathcal{A}$ often leads to slow convergence or even divergence. To address this problem, we implement a preconditioning strategy by transforming the systems into $\mathcal{P}^{-1}\mathcal{A}\mathbf{u}=\mathcal{P}^{-1}\mathbf{b}$ prior to applying Krylov subspace methods. Through appropriate parameter selection in the preconditioner $\mathcal{P}$, the transformed matrix $\mathcal{P}^{-1}\mathcal{A}$ achieves favorable eigenvalue clustering. This the enhancement of the spectral property enables Krylov subspace methods to attain significantly faster convergence rates when the preconditioned systems is solved.
\end{remark}

\section{Numerical experiments}

In this section, two numerical examples are given to demonstrate the  effectiveness of the preconditioner $\mathcal{P}$ over the existing preconditioners, i.e., $BS_2$ and BUT preconditioners, whose forms are given in (\ref{1.5}) and (\ref{BUT}). In order to compare them, we iteratively calculate the three-by-three block linear systems (\ref{1.4}) by the preconditioned GMRES methods under $BS_2$ and BUT preconditioners so as to obtain solution of the ILS problem (\ref{1.1}); similarly, we obtain the solution of the ILS problem (\ref{1.1}) by iteratively solving the block  three-by-three linear systems (\ref{1.3}) using the $\mathcal{P}$ preconditioning method. In the solution process, the initial vector $\mathbf{u}^{0}$ is chosen as $(\delta_1^{0};x^{(0)};\delta_2^{(0)})$ with $x^{(0)}=0$ and $(\delta_1^{(0)};\delta_2^{(0)})=b-Ax^{(0)}$,  the right-hand side vector $\mathbf{b}$ is chosen as the vector generated by free elements.  The iterations stop once the relative residual satisfies
\[
\mbox{RES}=\frac{\|\mathbf{b}-\mathcal{A}\mathbf{u}^{k}\|_2}{\|\mathbf{b}\|_2} < 10^{-8}
\]
with $\mathbf{u}^{k}$ the current approximation solution or if the prescribed maximum iteration count $k_{max}=1500$ is exceeded. The numerical behavior of the GMRES method are tested and evaluated in terms of the number of iteration steps (denoted `IT'), the computing times (denoted `CPU') and relative residual `RES'. The left preconditioned GMRES method (use the matlab function `\textbf{gmres}') is used to solve the block three-by-three linear systems (\ref{1.3}) and (\ref{1.4}).

We can find from the iteration process that the iteration solution $\mathbf{u}^{k}$ of the preconditioning methods contains the approximate solution $x^{k}$ of the ILS problem (\ref{1.1}). In order to reflect the effectiveness of the preconditioning method, the relative errors of the computed approximations $x^{k}$, denoted as `ERR', i.e.,
\[
\mbox{ERR}=\frac{\|x^{k}-x^{*}\|_{2}}{\|x^{*}\|_{2}}
\]
are reported with $x^{*}$ the initially exact solution.
All experiments are performed in MATLAB 2017(b) on an Intel Core(4G RAM) Windows 10 system.\\
\begin{example}\cite{Xin2025}
Our matrices $A_1$ are derived from the Tolosa Matrix collected by the Matrix Market, taken as TOLS340, TOLSs1090, TOLS2000, TOLS4000, WELL1033 and WELL1850. The first four Tolosa matrices arise in the stability analysis of a model of an airplane in flight, the problem has been analyzed at CERFACS in cooperation with the Aerospatiale Aircraft division. The last two matrices generate from the least squares problems in surveying. To ensure that our ILS problem has a unique solution, we let $A_2=0.3*eye(q,n)$, $b_1=rand(p,1)$, $b_2=rand(q,1)$. Then the ILS problem (\ref{1.1}) and the linear systems (\ref{1.3}) will be obtained.
\end{example}

\begin{figure}[!ht]
	\centering
	\includegraphics[height=4cm,width=5cm]{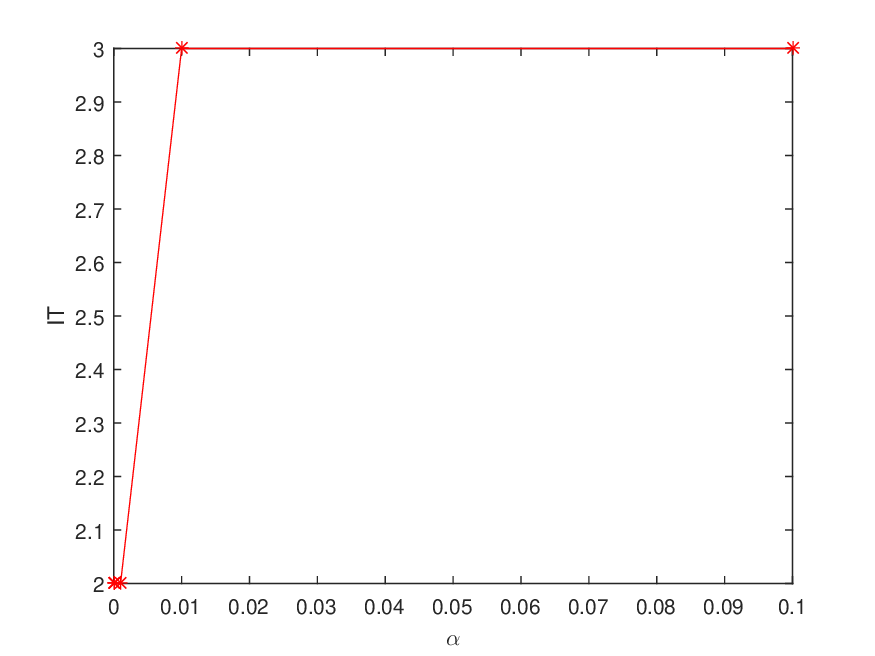}
	\includegraphics[height=4cm,width=5cm]{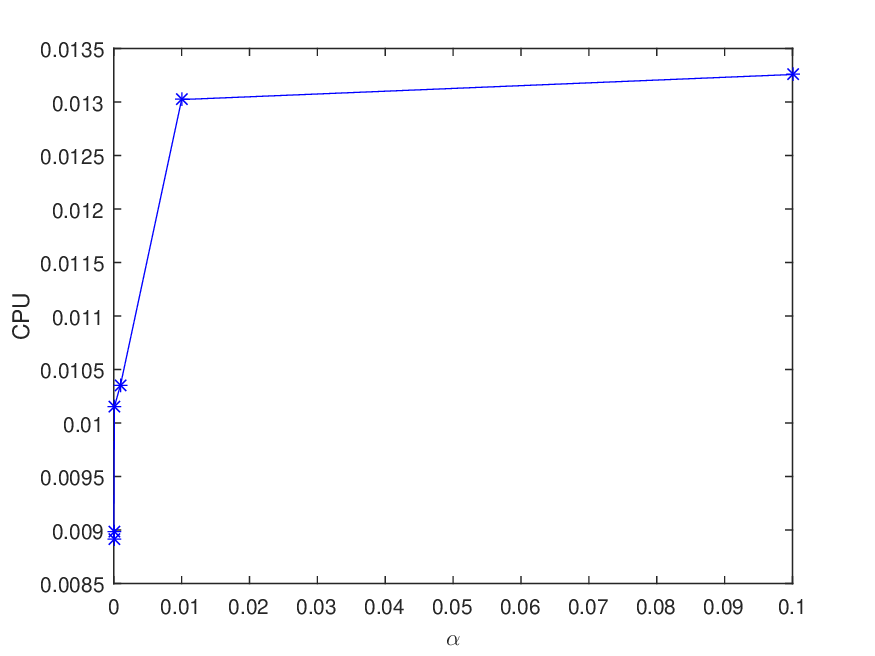}
	\caption{\footnotesize \emph{$\alpha$ vs IT and CPU of $\mathcal{P}$ preconditioned GMRES method under TOLS340 in Example 4.1.}}\label{fig1}
\end{figure}

For proposed preconditioner $\mathcal{P}$, we know from the theory that the parameter $\alpha$ should be selected close to $0$ such that the preconditioner is closer to the coefficient matrix of (\ref{1.3}). Therefore, in Figure 1, we plot the relationships between $\alpha$ and IT, CPU. It is obvious that the IT and CPU of the $\mathcal{P}$ preconditioned GMRES method are direct ratio with the value of the parameter. In order to better show the efficiency of $\mathcal{P}$, we can select $\alpha=10^{-6}$ in this example. Table \ref{tab1} shows that the no preconditioned GMRES method converges slowly; however, other preconditioning methods can converge quickly. Specifically speaking, $\mathcal{P}$ preconditioning method is less than the other two preconditioning methods in both iteration steps and CPU time. And the values of ERR show that the iteration solution obtained by $\mathcal{P}$ preconditioning method is closer to the exact solution. The above facts not only reflect the effectiveness of the selected parameter but also show the robustness of the preconditioner $\mathcal{P}$.

\begin{table}[!ht]
	\renewcommand\arraystretch{1}
	\centering
	\caption{Numerical results of the different preconditioned GMRES methods for Example 4.1.}\label{tab1}
 %with $p=n$, $q=10000$.}
{
\begin{tabular}{lccccc|cccccccccccccccccccc} \hline
&&&TOLS340&&&&&TOLS1090&&\\
&Pre.&I&$BS_2$&BUT&$\mathcal{P}$&Pre.&I&$BS_2$&BUT&$\mathcal{P}$\\\hline
&IT&90&4&4&2&IT&239&4&4&2\\
~&RES&6.0e-10&7.8e-10&2.0e-10&4.9e-14~&RES&1.3e-09&9.9e-10&2.1e-10&8.3e-16\\
~&ERR&1.1e-06&8.3e-10&8.3e-10&5.3e-14~&ERR&2.1e-06&5.0e-10&5.0e-10&9.6e-14\\
~&CPU&0.3247&0.0156&0.0147&{\bf 0.0097}~~&CPU&1.6779&0.0186&0.0171&{\bf0.0157}\\\hline
&&&TOLS2000&&&&&TOLS4000&&\\
&IT&419&4&4&2&IT&816&4&4&2\\
~&RES&6.7e-09&1.4e-09&3.2e-10&1.2e-15~&RES&6.1e-09&1.3e-09&9.5e-11&1.2e-13\\
~&ERR&6.5e-03&6.1e-10&6.1e-10&1.6e-13~&ERR&7.8e-02&7.2e-10&5.3e-10&3.7e-13\\
~&CPU&5.6294&0.0237&0.0234&{\bf0.0229}&~CPU&25.9522&0.0698&0.0609&{\bf0.0482}\\\hline
%&&&ILLC1033&&&&&ILLC1850&&\\
%&IT&117&112&112&2~&IT&919&229&229&3\\
%~&RES&4.1e-09&1.7e-08&6.7e-09&2.0e-09~&RES&9.6e-09&8.4e-09&7.8e-09&5.7e-11\\
%~&ERR&3.7e-08&4.4e-03&2.4e-03&3.3e-09~&ERR&5.4e-08&2.2e-05&2.0e-05&5.7e-11\\
%~&CPU&0.4531&0.7051&0.6725&{\bf0.0218}~&CPU&26.6778&2.8525&2.7749&{\bf0.0532}\\\hline
&&&WELL1033&&&&&WELL1850&&\\
&IT&124&40&40&2&IT&621&133&133&3\\
~&RES&5.6e-09&2.8e-09&1.2e-09&3.3e-09~&RES&9.4e-09&3.6e-09&7.8e-09&1.2e-11\\
~&ERR&1.8e-08&2.7e-08&2.2e-08&4.6e-09~&ERR&2.4e-07&1.0e-07&3.0e-07&1.3e-11\\
~&CPU&0.4948&0.1631&0.1500&{\bf0.0203}~&CPU&11.7769&1.1656&1.0599&{\bf0.0486}\\\hline
\end{tabular}}

{\caption{The different numerical results of the preconditioned GMRES methods for Example 4.2.}\label{tab2}
\begin{tabular}{l|cccc|cccc}
\hline
&\multicolumn{4}{c|}{$p=256,q=100,n=128$}&\multicolumn{4}{c}{$p=300,q=212,n=256$}\\
\hline
Pre.&IT &RES&ERR&CPU &IT &RES&ERR&CPU \\
\hline
No&51&9.2e-09&1.3e-04&0.0254&76&9.8e-09&8.6e-04&0.0891\\
$BS_2$&3&7.8e-09&3.5e-09&0.0119&5&1.5e-10&1.8e-10&0.0625\\
$\mathrm{BUT}$&3&1.2e-09&9.9e-10&0.0113&4&8.8e-09&8.9e-09&0.0483\\
$\mathcal{P}$&2&3.4e-14&2.0e-12&{\bf0.0097}&2&2.2e-12&1.2e-11&{\bf0.0404}\\
\hline
&\multicolumn{4}{|c}{$p=524,q=500,n=512$}&\multicolumn{4}{|c}{$p=1048,q=1000,n=1024$}\\
\hline
Pre&IT &RES&ERR&CPU &IT &RES&ERR&CPU \\
\hline
No&158&9.3e-09&1.9e-02&0.3868&325&9.9e-09&9.3e-01&4.1752\\
$BS_2$&6&5.4e-10&5.5e-10&0.2533&12&8.5e-09&9.4e-09&2.2492\\
$\mathrm{BUT}$&6&6.3e-10&8.3e-10&0.2361&12&6.9e-09&4.7e-09&2.1866\\
$\mathcal{P}$ &2&9.4e-09&9.2e-09&{\bf0.1093}& 5&3.7e-10&6.7e-10&{\bf1.0932}\\
\hline
\end{tabular}}
\end{table}

\begin{example}\cite{Liu2011,Xin2025}
Let $p,~q,~ n$ take different values, $\epsilon=10^{-4}$ and
\begin{align*}
	\tilde{B}=Y\left(\begin{matrix}
		D \\
		0\\
	\end{matrix}\nonumber
	\right)Z^{T}\in  \mathbb{R}^{p\times n},
\end{align*}
where $Y\in  \mathbb{R}^{p\times p}$, $Z\in  \mathbb{R}^{n\times n}$ are given orthogonal matrices and $ D={\rm diag}(1,\frac{1}{2},..., \frac{1}{n})\in  \mathbb{R}^{n\times n}$. Next let $B=\tilde{B}+\epsilon E$, $d=\tilde{B}1_{n}+\epsilon f$, where $1_{n}$ denotes an $n\times 1$ column vector of all ones, $E$ and $f$ are given error matrix and vector generated by using the Matlab.

If $B^{T}B-\sigma_{n+1}^{2}I_{n}$ is positive definite, the solution of the total least squares problem \cite{Bjorck1996} related to $B$, $d$ is as follows
\begin{equation*}
	x_{TLS}=(B^{T}B-\sigma_{n+1}^{2}I_{n})^{-1}B^{T}d,
\end{equation*}
where $\sigma_{n+1}$ is the smallest singular value of $(B, d)$. It is equivalent to solving the ILS problem with
\begin{align*}
A=\left(\begin{matrix}
	B \\
	\sigma_{n+1}I\\
\end{matrix}\nonumber
\right),\quad
b=\left(\begin{matrix}
	d \\
	0\\
\end{matrix}\nonumber
\right).
\end{align*}
\end{example}

%\begin{table}[!ht]
%\renewcommand\arraystretch{1}
%{\small
%\caption{The different numerical results of the preconditioned GMRES methods for example 4.2 with $tol=10^{-6}$.}\label{tab2}
%\begin{center}
%{
%\begin{tabular}{l|cccc|cccc}
%\hline
%&\multicolumn{4}{c|}{$p=256,q=100,n=128$}&\multicolumn{4}{c}{$p=300,q=212,n=256$}\\
%\hline
%Pre.&IT &RES&ERR&CPU &IT &RES&ERR&CPU \\
%\hline
%No&51&9.2e-09&1.3e-04&0.0254&76&9.8e-09&8.6e-04&0.0891\\
%$BS_2$&3&7.8e-09&3.5e-09&0.0119&5&1.5e-10&1.8e-10&0.0625\\
%$\mathrm{BUT}$&3&1.2e-09&9.9e-10&0.0113&4&8.8e-09&8.9e-09&0.0483\\
%$\mathcal{P}$&2&3.4e-14&2.0e-12&{\bf0.0097}&2&2.2e-12&1.2e-11&{\bf0.0404}\\
%\hline
%&\multicolumn{4}{|c}{$p=524,q=500,n=512$}&\multicolumn{4}{|c}{$p=1048,q=1000,n=1024$}\\
%\hline
%Pre&IT &RES&ERR&CPU &IT &RES&ERR&CPU \\
%\hline
%No&158&9.3e-09&1.9e-02&0.3868&325&9.9e-09&9.3e-01&4.1752\\
%$BS_2$&6&5.4e-10&5.5e-10&0.2533&12&8.5e-09&9.4e-09&2.2492\\
%$\mathrm{BUT}$&6&6.3e-10&8.3e-10&0.2361&12&6.9e-09&4.7e-09&2.1866\\
%$\mathcal{P}$ &2&9.4e-09&9.2e-09&{\bf0.1093}& 5&3.7e-10&6.7e-10&{\bf1.0932}\\
%\hline
%\end{tabular}}
%\end{center}}
%\end{table}

\begin{figure}[!ht]
\centering
	\subfigure[\label{fig:a}]{
		\includegraphics[height=3.5cm,width=5cm]{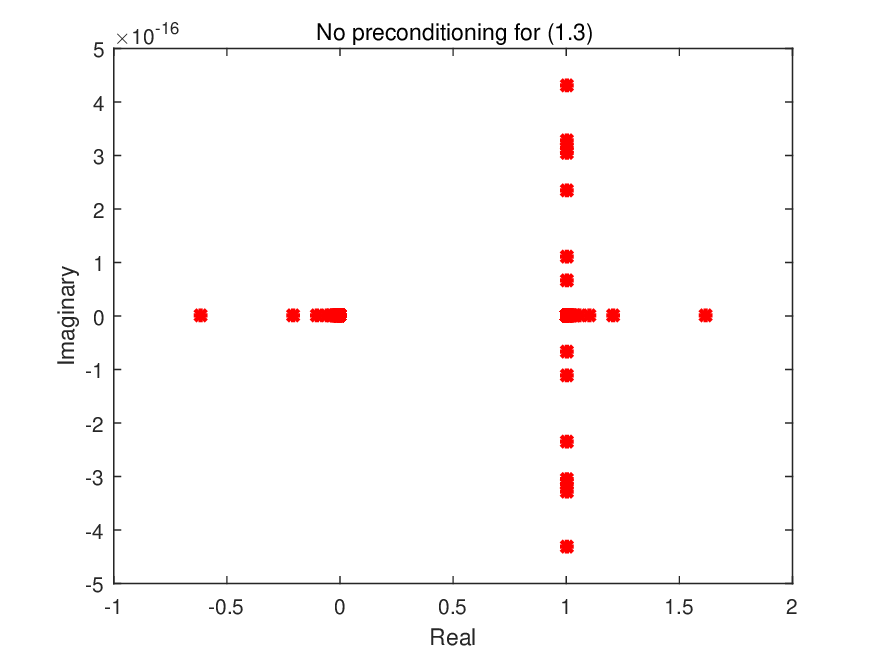}}
	\subfigure[\label{fig:b}]{
		\includegraphics[height=3.5cm,width=5cm]{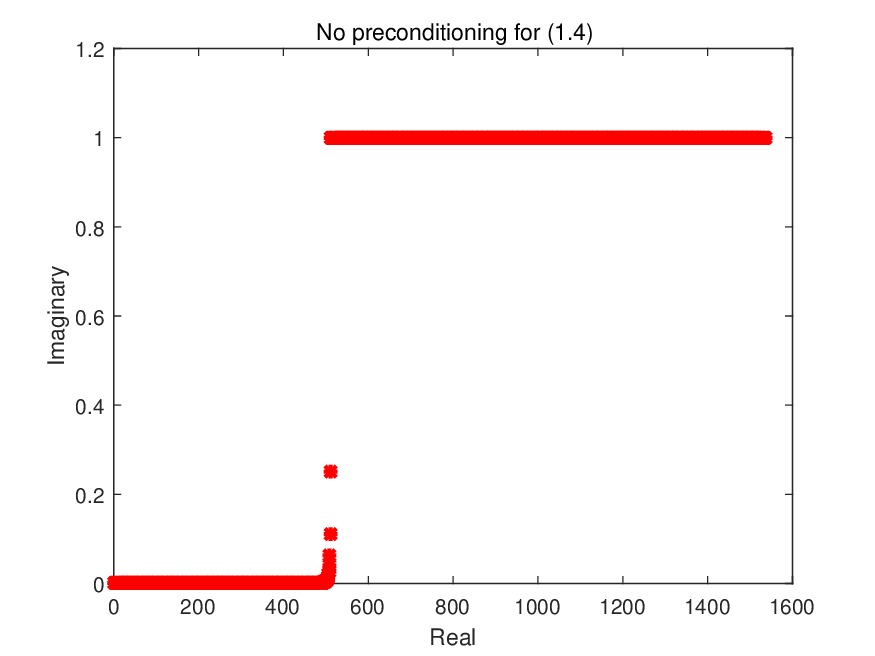}}\\
	\subfigure[\label{fig:c}]{
		\includegraphics[height=3.5cm,width=5cm]{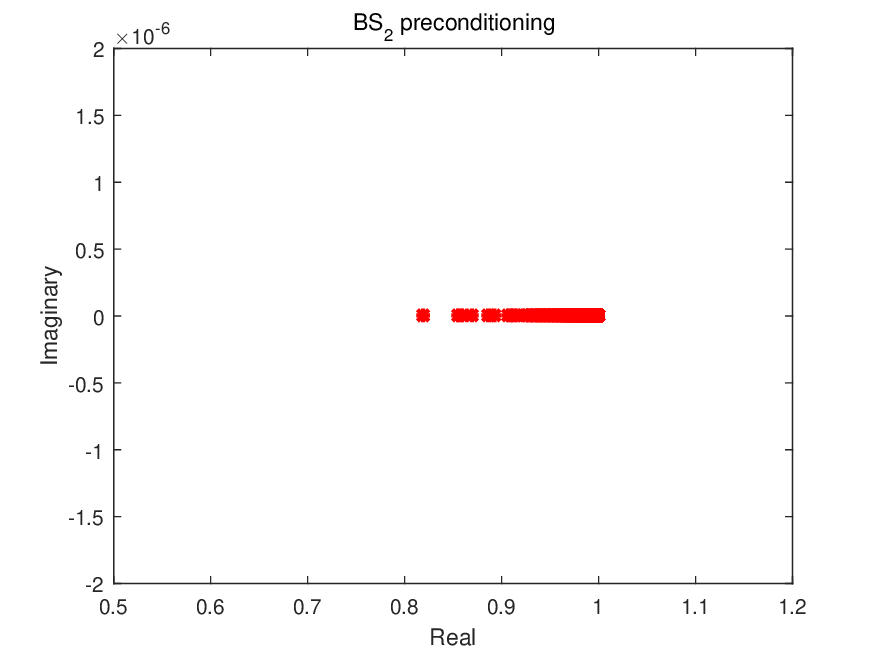}}
	\subfigure[\label{fig:d}]{
		\includegraphics[height=3.5cm,width=5cm]{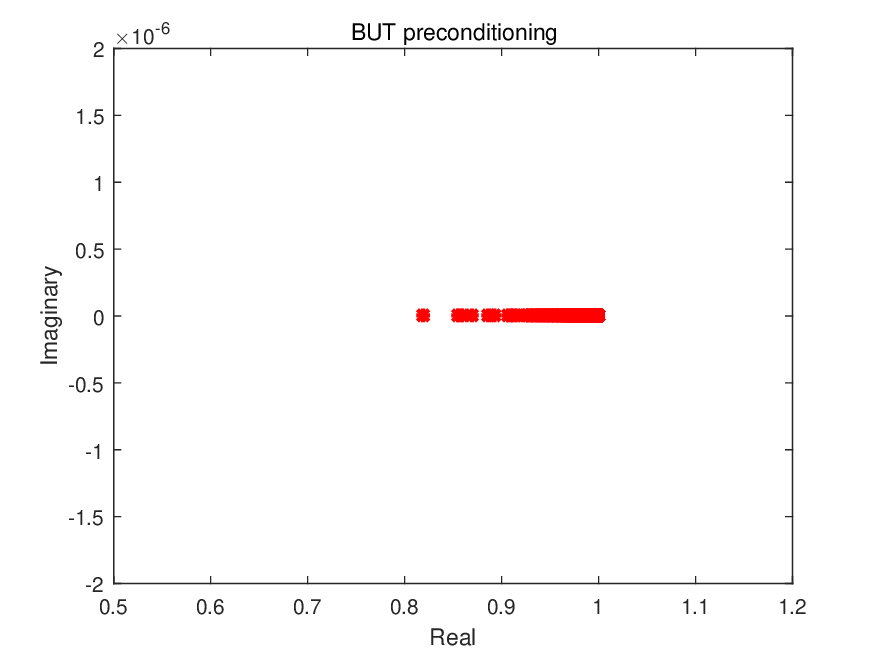}}
	\subfigure[\label{fig:e}]{
		\includegraphics[height=3.5cm,width=5cm]{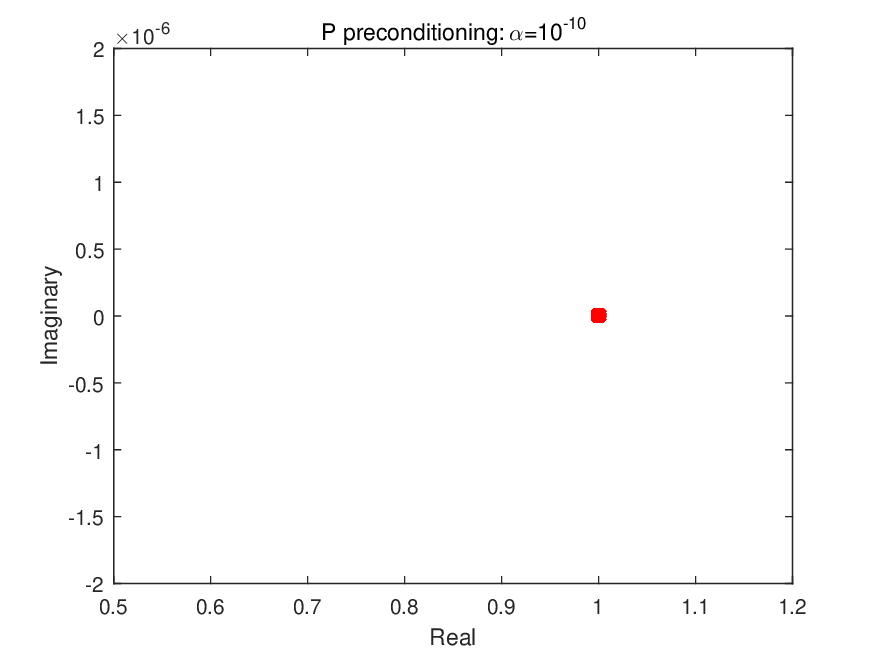}}
\caption{\footnotesize \emph{The eigenvalue distribution of the initial coefficient matrices $\mathcal{\hat{A}}$, $\mathcal{A}$ of (\ref{1.3}), (\ref{1.4}) and the preconditioned matrices under different preconditioners in example 4.2 with $p=524,q=500,n=512$. }}
\end{figure}

Similar to example 4.1, we can select $\alpha=10^{-10}$ in example 4.2 so as to have good convergence of the $\mathcal{P}$ preconditioned GMRES method.  We find from table 2 that the $\mathcal{P}$ preconditioned GMRES method is better than $BS_2$ and BUT preconditioned GMRES methods in both iteration steps and CPU time under various dimensions. And its iteration solution can also approximate the exact solution well, which is enough to prove the advantage of $\mathcal{P}$ preconditioning for the initial linear systems (\ref{1.3}). Besides, in figure 2, we plot the eigenvalues distribution of the initial coefficient matrix $\mathcal{A}$, $\mathcal{\hat{A}}$ and different preconditioned matrices in example 4.2. Figures a and b shows that the eigenvalues of the coefficient matrices $\mathcal{A}$ and $\mathcal{\hat{A}}$ are scattering, however, preconditioning could better improve the eigenvalues distribution of matrices. Besides, the eigenvalues of the $\mathcal{P}$ preconditioned matrix gather at the point $(1,0)$ under $\alpha=10^{-10}$, the tight spectrums of the preconditioned matrix lead to stable numerical performances, see tables for details.

%We compare GMRES method with the preconditioned GMRES method with preconditioner $\mathcal{P}$ from aspects of the number of iteration steps (denoted by \lq IT\rq) and elapsed CPU times in seconds (denoted by \lq CPU\rq).

\section{Conclusions}
In this paper, a new preconditioner was proposed to solve the  block three-by-three linear systems derived from the ILS problem, and its efficiency was more than the serval existing preconditioners in \cite{Li2024,Xin2025}. The future work will focus on the efficient algorithms for solving ILS problem derived from the practical problems.

\section*{Declarations}

%\noindent{\bf Ethical Approval} Not Applicable.\\

\noindent{\bf Availability of supporting data} Data availability is not applicable to this article as no new data were created or analyzed in this study.\\

%\noindent{\bf Competing interests} The authors declare no competing interests.\\

%\noindent{\bf Authors' contributions} Li and Meng wrote the main manuscript text, Xin tested the numerical experiments.\\

\noindent{\bf Funding}
This work was supported by the National Natural Science Foundation of China (No.12361082).

\end{document}